\documentclass[12pt,reqno]{amsart}

\usepackage[arrow,matrix,curve]{xy}

\usepackage[dvips]{graphicx} 

\usepackage{amssymb, latexsym, amsmath, amscd, array,
}

\theoremstyle{definition}

\numberwithin{equation}{section}

\newcommand\N {{\mathbb N}} 

\newcommand\R {{\mathbb R}}
\newcommand\Q {{\mathbb Q}}

\newcommand\RRR{\mbox{I\!I\!R}}

\newcommand\RRRhuge{{\Huge\mbox{I\!I\!R}}}

\newcommand\Los{{\L}o{\'s}}

\begin{document}

\thispagestyle{empty}

\title[Who gave you the Cauchy--Weierstrass tale?]{Who gave you the
Cauchy--Weierstrass tale?  The dual history of rigorous calculus}

\author{Alexandre Borovik}

\address{School of Mathematics, University of Manchester, Oxford
Street, Manchester, M13 9PL, United Kingdom}
\email{alexandre.borovik@manchester.ac.uk}

\author{Mikhail G. Katz$^{0}$}

\address{Department of Mathematics, Bar Ilan University, Ramat Gan
52900 Israel} \email{katzmik@macs.biu.ac.il}

\footnotetext{Supported by the Israel Science Foundation grant
1294/06.}



\begin{abstract}
Cauchy's contribution to the foundations of analysis is often viewed
through the lens of developments that occurred some decades later,
namely the formalisation of analysis on the basis of the epsilon-delta
doctrine in the context of an Archimedean continuum.  What does one
see if one refrains from viewing Cauchy as if he had read Weierstrass
already?  One sees, with Felix Klein, a parallel thread for the
development of analysis, in the context of an infinitesimal-enriched
continuum.  One sees, with Emile Borel, the seeds of the theory of
rates of growth of functions as developed by Paul du Bois-Reymond.
One sees, with E.~G.~Bj\"orling, an infinitesimal definition of the
criterion of uniform convergence.  Cauchy's foundational stance is
hereby reconsidered.
\end{abstract}

\maketitle

\tableofcontents

\section{Introduction}
\label{one}

When the second-named author first came across a recent book by
Anderson {\em et al\/} entitled {\em Who Gave You the Epsilon?  And
Other Tales of Mathematical History\/}~\cite{AKW}, he momentarily
entertained a faint glimmer of hope.  The book draws its title from an
older essay, entitled {\em Who Gave You the Epsilon?  Cauchy and the
origins of rigorous calculus\/}~\cite{Grab}.  The faint hope was that
the book would approach the thesis implied by the title of the older
essay, in a critical spirit, namely, as an ahistorical%
\footnote{See Grattan-Guinness's comment on ahistory in the main text
around footnote~\ref{muddle}.}
{\em tale\/} in need of re-examination.  Anderson~{\em et al\/} not
having undertaken the latter task, such an attempt is made here.

Cauchy based his definitions of both limits and infinitesimals on the
concept of a variable quantity.%
\footnote{\label{gil}C.~Gilain's claims to the effect that ``Cauchy
d\'efinissait le concept d'infiniment petit \`a l'aide du concept de
limite, qui avait le premier r\^ole'' \cite[footnote~67]{Gil} are both
in error; see Subsection~\ref{gilain} below.  Already in 1973, Hourya
Benis Sinaceur wrote: ``on dit trop rapidement que c'est Cauchy qui a
introduit la `m\'ethode des limites' entendant par l\`a, plus ou moins
vaguement, l'emploi syst\'ematique de l'{\em \'epsilonisation\/}''
\cite[p.~108]{Si}, and pointed out that Cauchy's definition of limit
resembles, not that of Weierstrass, but rather that of Lacroix dating
from 1810 \cite[p.~109]{Si}.}
The variable quantities used in his {\em Cours d'Analyse\/} in 1821
are generally understood by scholars as being (discrete) sequences of
values.  Cauchy wrote that a variable quantity tending to zero becomes
infinitely small.  How do his null sequences become infinitesimals?
In 1829, Cauchy developed a detailed theory of infinitesimals of
arbitrary order (not necessarily integer), based on rates of growth of
functions.  How does his theory connect with the work of later
authors?  In 1902, E.~Borel elaborated on du Bois-Reymond's theory of
rates of growth, and outlined a general ``theory of increase'' of
functions, as a way of implementing an infinitesimal-enriched
continuum.  Borel traced the lineage of such ideas to Cauchy's text.
We examine several views of Cauchy's foundational contribution in
analysis.

Cauchy's foundational stance has been the subject of an ongoing
controversy.  The continuing relevance of Cauchy's foundational stance
stems from the fact that Cauchy developed some surprisingly modern
mathematics using infinitesimals.  The following four items deserve to
be mentioned:

(1) {\em Cauchy's proof of the binomial formula (series) for arbitrary
  exponents\/}.  The proof exploits infinitesimals.  D.~Laugwitz
  \cite[p.~266]{Lau87} argues that this is the first correct proof of
  the formula.%
\footnote{except for Bolzano's proof in 1816, see
\cite[p.~657]{Lau97}.}

(2) {\em Cauchy's use of the Dirac delta function\/}.  Over a century
before Dirac, Cauchy exploited ``delta functions'' to solve problems
in Fourier analysis and in the evaluation of singular integrals.  Such
functions were defined in terms of an infinitesimal parameter, see
Freudenthal~\cite{Fr71a}, Laugwitz \cite[p.~219]{Lau89}, \cite{Lau92}.

(3) {\em Cauchy's definition of continuity\/}.  Cauchy defined
continuity of a function~$y=f(x)$ as follows: {\em an infinitesimal
  change~$\alpha$ of the independent variable~$x$ always produces an
  infinitesimal change
\hbox{$f(x+\alpha)-f(x)$} of the dependent variable~$y$\/}
\cite[p.~34]{Ca21}.  Nearly half a century later, Weierstrass
reconstructed Cauchy's infinitesimal definition in the following
terms: for every~$\epsilon>0$ there exists a~$\delta>0$ such that for
every real~$\alpha$, if~$|\alpha|<\delta$
then~$|f(x+\alpha)-f(x)|<\epsilon$.  Many historians have sought to
interpret Cauchy's definition as a proto-Weierstrassian definition of
continuity in terms of limits.%
\footnote{\label{smi}Thus, Smithies \cite[p.~53, footnote~20]{Sm}
cites the {\em page\/} in Cauchy's book where Cauchy gave the
infinitesimal definition, but goes on to claim that the concept of
{\em limit\/} was Cauchy's ``essential basis'' for his concept of
continuity \cite[p.~58]{Sm}.  Smithies looked in Cauchy, saw the
infinitesimal definition, and went on to write in his paper that he
saw limits (see Subsection~\ref{gilain} for an analysis of a similar
misconception in Gilain).  Such automated translation has been common
at least since Boyer \cite[p.~277]{Boy}.}

(4) {\em Cauchy's ``sum theorem".\/} This result asserts the
convergence to a {\em continuous\/} function, of a series of
continuous functions under a suitable condition of convergence.  The
sum theorem has been the subject of a historical controversy, ever
since A.~Robinson \cite[p.~271-273]{Ro66} proposed a novel reading of
the sum theorem that would make it correct.  The controversy hinges on
the question whether the convergence condition was meant by Cauchy to
hold at the points of an Archimedean continuum, or at the points of a
Bernoullian continuum, namely, an infinitesimal-enriched continuum
(see Subsection~\ref{micro1}).  Lakatos presented a paper%
\footnote{The paper was published posthumously in 1978 and edited by
J.~Cleave \cite{La78}.}
in 1966 where he argues that the 1821 result is correct as stated.%
\footnote{Lakatos thereby reversed his position as presented in his
{\em Proofs and Refutations\/}~\cite{La76}.}
Laugwitz concurs, see e.g. his 1989 text \cite{Lau89}.
Post-Weierstrassian historians tend both (a) to reject Cauchy's
infinitesimals, claiming they are merely shorthand for limits, and (b)
to claim that his 1821 sum theorem was false.  The 1821 formulation of
the sum theorem may in the end be too ambiguous to know what Cauchy's
intention was at the time, if Cauchy himself knew.%
\footnote{See more details in Subsection~\ref{micro1}.}

The four examples given above hopefully illustrate the interest in
understanding the nature of Cauchy's foundational contribution.

In Section~\ref{stance}, we provide a re-appraisal of Cauchy's
foundational stance, including Cauchy's theory of orders of
infinitesimals based on orders of growth of functions
(Subsection~\ref{orders}); how Cauchy null sequences become
infinitesimals (Subsection~\ref{second}); a detailed textual study of
the related portions of the {\em Cours d'Analyse\/}
(Subsection~\ref{two}); an analysis of Cauchy's key term {\em
toujours\/}, strengthening the hypothesis of his sum theorem
(Subsection~\ref{micro1}); an examination of some common
misconceptions in the Cauchy literature (Subsection~\ref{25}).

In Section~\ref{modern}, we analyze the views of some modern authors.
In Section~\ref{time}, we provide a timeline of the development of
theories of infinitesimals based on orders of growth of functions,
from Cauchy through Paul du Bois-Reymond to modern times.  Some
conclusions and outlook for the future may be found in
Section~\ref{five}.  Appendix~\ref{rival} outlines the relevant
mathematical material on the rival continua, Archimedean and
Bernoullian.

\section{A reappraisal of Cauchy's foundational stance}
\label{stance}

\subsection{Cauchy's theory of orders of infinitesimals}
\label{orders}

Cauchy's 1821 {\em Cours d'Analyse\/} \cite{Ca21} presented only a
theory of infinitesimals of polynomial rate of growth as compared to a
given ``base'' infinitesimal~$\alpha$.  

The shortcoming of such a theory is its limited flexibility.  Since
Cauchy only considers infinitesimals behaving as polynomials of a
fixed infinitesimal, called the ``base" infinitesimal in 1823
\cite{Ca23}, his framework imposes obvious limitations on what can be
done with such infinitesimals.  Thus, one typically can't extract the
square root of such a ``polynomial'' infinitesimal.

What is remarkable is that Cauchy did develop a theory to overcome
this shortcoming.  Cauchy's theory of infinitesimals of arbitrary
order (not necessarily integer) is noted by A.~Borel
\cite[p.~35-36]{Bo02} (see below).

In 1823, and in more detail in 1829, Cauchy develops a more flexible
theory, where an infinitesimal is represented by an arbitrary {\em
function\/} (rather than merely a polynomial) of a base infinitesimal,
denoted ``$i$", see Chapter 6 in \cite{Ca29}.  The title of the
chapter is significant.  Indeed, the title refers to functions as {\em
representing\/} the infinitesimals; more precisely, {\em ``fonctions
qui repr\'esentent des quantit\'es infiniment petites"\/}.  Here is
what Cauchy has to say in~1829:
\begin{quote}
Designons par~$a$ un nombre constant, rationnel ou irrationnel; par
$i$ une quantit\'e infiniment petite, et par~$r$ un nombre variable.
Dans le systeme de quantit\'es infiniment petites dont~$i$ sera la
base, une fonction de~$i$ represent\'ee par~$f(i)$ sera un infiniment
petit de l'ordre~$a$, si la limite du rapport~$f(i)/i^r$ est nulle
pour toutes les valeurs de~$r$ plus petites que~$a$, et infinie pour
toutes les valeurs de~$r$ plus grandes que~$a$ \cite[p. 281]{Ca29}.
\end{quote}

Laugwitz \cite[p.~271]{Lau87} explains this to mean that the order~$a$
of the infinitesimal~$f(i)$ is the uniquely determined real number
(possibly~$+\infty$, as with the function~$e^{-1/t^2}$) such
that~$f(i)/i^r$ is infinitesimal for~$r < a$ and infinitely large
for~$r > a$.

Laugwitz \cite[p.~272]{Lau87} notes that Cauchy provides an example of
functions defined on positive reals that represent infinitesimals of
orders~$\infty$ and~$0$, namely 
\[
e^{-1/i} \mbox{ \ and \ } \frac{1}{\log i}
\]
(see Cauchy \cite[p.~326-327]{Ca29}).

The development of non-Archimedean systems based on orders of growth
was pursued in earnest at the end of the 19th century by such authors
as Stolz, du Bois-Reymond, Levi-Civita, and Borel.%
\footnote{\label{ehr}See P.~Ehrlich's \cite{Eh06} for a more detailed
discussion, as well as Section~\ref{time} below.}
Such systems have an antecedent in Cauchy's theory of infinitesimals
as developed in his texts dating from 1823 and~1829.  Thus, in 1902,
E.~Borel \cite[p.~35-36]{Bo02} elaborated on du Bois-Reymond's theory
of rates of growth, and outlined a general ``theory of increase'' of
functions, as a way of implementing an infinitesimal-enriched
continuum.  Borel traced the lineage of such ideas to an 1829 text of
Cauchy's on the rates of growth of functions, see Fisher
\cite[p.~144]{Fi} for details.  In 1966, A.~Robinson pointed out that
\begin{quote}
Following Cauchy's idea that an infinitely small or infinitely large
quantity is associated with the behavior of a function~$f(x)$, as~$x$
tends to a finite value or to infinity, du Bois-Reymond produced an
elaborate theory of orders of magnitude for the asymptotic behavior of
functions \dots Stolz tried to develop also a theory of arithmetical
operations for such entities \cite[p.~277-278]{Ro66}.
\end{quote}
Robinson traces the chain of influences further, in the following
terms:
\begin{quote}
It seems likely that Skolem's idea to represent infinitely large
natural numbers by number-theoretic functions which tend to infinity
(Skolem [1934]),%
\footnote{\label{1966a}The reference is to Skolem's 1934 work
\cite{Sk}.  The evolution of modern infinitesimals is traced further
in the main text in Section~\ref{time} following footnote~\ref{1966}.}
also is related to the earlier ideas of Cauchy and du Bois-Reymond
\cite[p.~278]{Ro66}.
\end{quote}

Cauchy's approach is by no means a variant of Robinson's approach.%
\footnote{More specifically, Cauchy was not in the possession of the
mathematical tools required to either formulate or justify the
ultrapower construction, requiring as it does a set-theoretic
framework (dating from the end of the 19th century) together with the
existence of ultrafilters (not proved until 1930 by Tarski
\cite{Tar}), see Appendix~\ref{rival}.}
Similarly, the notion of Cauchy as a pre-Weierstrassian who allegedly
``tried to dislodge infinitesimals from analysis'' is just as dubious.
Taken to its logical conclusion, the dogma of Cauchy as a
pre-Weierstrassian can assume comical proportions.  Thus, the
fashionable Stephen Hawking comments that Cauchy
\begin{quote}
was particularly concerned to banish infinitesimals \cite[p.~639]{Ha},
\end{quote}
yet {\em on the very same page\/} 639, Hawking quotes Cauchy's {\em
infinitesimal\/} definition of continuity in the following terms:
\begin{quote}
the function~$f(x)$ remains continuous with respect to~$x$ between the
given bounds, if, between these bounds, an infinitely small increment
in the variable always produces an infinitely small increment in the
function itself \cite[p.~639]{Ha}.
\end{quote}
Did Cauchy ``banish'' infinitesimals?  Using infinitely small
increments is an odd way of doing so.  Similarly, historian J.~Gray
lists {\em continuity\/} among concepts Cauchy allegedly defined
\begin{quote}
using careful, if not altogether unambiguous, {\bf limiting} arguments
\cite[p.~62]{Gray08} [emphasis added--authors],
\end{quote}
whereas in point of fact, {\em limits\/} appear in Cauchy's definition
only in the sense of the {\em endpoints\/} of the domain of
definition.  Analogous misconceptions in Gilain are analyzed in
Subsection~\ref{gilain}.%
\footnote{\label{gray}See \cite{KK11a, KK11b} for further discussion
of the post-Weierstrassian bias in Cauchy scholarship.}
Similar post-Weierstrassian sentiments were expressed by editor Michel
Blay's referee for the periodical {\em Revue d'histoire des
sciences\/} where the present article was submitted for publication in
2010, and rejected based on the following assessment by referee 1:%
\footnote{The pdf of the version submitted to Blay, as well as the two
referee reports, may be found at
http://u.cs.biu.ac.il/$\sim$katzmik/straw.html}
\begin{quote}
Our author interprets A. Cauchy's approach as a formation of the idea
of an infinitely small - a variant of the approach which was developed
in the XXth century in the framework of the nonstandard analysis (a
hyperreal version of E. Hewitt, J. Los, A. Robinson).
\end{quote}
The referee's summary is an inaccurate description of the text
refereed.  Cauchy's approach is not a variant of Robinson's, and was
never claimed to be in this text.  Based on such a strawman version of
the article's conception, the referee came to the following
conclusion:
\begin{quote}
From my point of view the author's arguments to support this
conception are quite unconvincing.
\end{quote}
Indeed, we find such such a strawman conception unconvincing, but the
conception was the referee's, not the authors'.  The referee concluded
as follows:
\begin{quote}
The fact that the actual infinitesimals lived somewhere in the
consciousness of A. Cauchy (as in many another mathematicians of XIXth
- XXth centuries as, for example, N.N. Luzin) does not abolish his
(and theirs) constant aspiration to dislodge them in the
subconsciousness and to found the calculus on the theory of limit.%
\footnote{Opposing infinitesimals and limits is in itself a conceptual
error; see footnote~\ref{ult}.}
\end{quote}
Did Cauchy have a ``constant aspiration to dislodge infinitesimals in
the subconsciousness and to found the calculus on the theory of
limit''?  Certainly not.  Felix Klein knew better: fifty years before
Robinson, Klein clearly realized the potential of the infinitesimal
approach to the foundations.  Having outlined the developments in real
analysis associated with Weierstrass and his followers, Klein pointed
out that
\begin{quote}
The scientific mathematics of today is built upon the series of
developments which we have been outlining.  But an essentially
different conception of infinitesimal calculus has been running
parallel with this [conception] through the centuries
\cite[p.~214]{Kl}.
\end{quote}
Such a different conception, according to Klein, ``harks back to old
metaphysical speculations concerning the structure of the continuum
according to which this was made up of [\dots] infinitely small parts''
\cite[p.~214]{Kl}.%
\footnote{\label{klein2}Klein also formulated a criterion of what it
would take for a theory of infinitesimals to be successful.  Namely,
one must be able to prove a mean value theorem for arbitrary
intervals, including infinitesimal ones.  In 1928, A.~Fraenkel
\cite[pp.~116-117]{Fran} formulated a similar requirement in terms of
the mean value theorem.  Such a Klein-Fraenkel criterion is satisfied
by the Hewitt-\Los-Robinson theory by the transfer principle, see
Appendix~\ref{rival}.}

Cauchy did not aspire to dislodge infinitesimals from analysis; on the
contrary, he used them with increasing frequency in his work,
including his 1853 article \cite{Ca53} where he relied on
infinitesimals to express the property of uniform convergence, as we
analyze in Subsection~\ref{micro1}.

\subsection{How does a null sequence become a Cauchy infinitesimal?}
\label{second}

The nature of Cauchy's infinitesimals has been the subject of an
ongoing debate for a number of decades.  Post-Weierstrassian
historians tend to dismiss Cauchy's {\em infiniment petits\/} as
merely a linguistic device masking Cauchy's use of the limit concept.
From this perspective, Cauchy's infinitesimals are alleged to
represent an early anticipation of the more rigorous methods developed
in the second half of the 19th century, namely epsilontic analysis.
Hourya Benis Sinaceur~\cite{Si} presented a critical analysis of the
proto-Weierstrassian approach to Cauchy in 1973, pointing out in
particular that Cauchy's notion of limit is for the most part a
kinetic one rather than an epsilontic one.%
\footnote{\label{pourciau}Pourciau \cite{Pou} argues that Newton
possessed a clear kinetic conception of limit similar to Cauchy's, and
cites Newton's lucid statement to the effect that ``Those ultimate
ratios \dots are not actually ratios of ultimate quantities, but
limits \dots which they can approach so closely that their difference
is less than any given quantity\dots'' See Newton \cite[p.~39]{New46}
and \cite[p.~442]{New99}.  The same point, and the same passage from
Newton, appeared a century earlier in Russell \cite[item 316,
p.~338-339]{Ru03}.  See also footnote~\ref{ult}.}
Other historians have taken Cauchy's infinitesimals at face value, see
e.g., Lakatos \cite{La78}, Laugwitz~\cite{Lau89}.  The studies in the
past decade include Sad, Teixeira, and Baldino~\cite{STB}, Br\aa ting
\cite{Br07}, and Katz \& Katz \cite{KK11b}.

To what extent did Cauchy intend the process that he described as a
null sequence ``becoming" an infinitesimal, to involve some kind of a
collapsing?

Cauchy did not have access to the modern set theoretic mentality
(currently dominant in the area of mathematics), where equivalence
relation and quotient space constructions are taken for granted.  One
can still ponder the following question: to what extent may Cauchy
have anticipated such collapsing phenomena?

Part of the difficulty in answering such a question is Cauchy's
hands-on approach to foundations.  Cauchy was less interested in
foundational issues than, say, Bolzano.%
\footnote{Thus, Hourya Benis Sinaceur writes: ``Bolzano est avant tout
pr\'eoccup\'e de rigueur th\'eorique, d'o\`u le soucit constant de
d\'emonstrations formelles [\dots] Cauchy, au contraire, donne
rarement des d\'emonstrations en forme [\ldots] Son expos\'e a
davantage des qualit\'es de {\em synth\`ese\/} que de rigueur
formelle'' \cite[p.~102]{Si}.}
To Cauchy, getting your sedan out of the garage was the only
justification for shoveling away the snow that blocks the garage door.
Once the door open, Cauchy is in top gear within seconds, solving
problems and producing results.  To him, infinitesimals were the
asphalt under the snow, not the snow itself.  Bolzano wanted to shovel
away {\em all\/} the snow from the road.  Half a century later, the
triumvirate%
\footnote{Boyer \cite[p.~298]{Boy} refers to Cantor, Dedekind, and
Weierstrass as ``the great triumvirate''.}
shovel ripped out the asphalt together with the snow, intent on
consigning the infinitesimal to the dustbin of history.

For these reasons, it is not easy to gauge Cauchy's foundational
stance precisely.  For instance, is there evidence that he felt that
two null sequences that coincide except for a finite number of terms,
would ``generate" the same infinitesimal?%
\footnote{Note that the term ``generate'' was used by K.~Br\aa
ting~\cite{Br07} to describe the passage from sequence to infinitesimal
in Cauchy.}
Support for this idea comes from several sources, including a very
unlikely one, namely Felscher's essay~\cite{Fe} in {\em The American
Mathematical Monthly\/} from 2000, where he attacks both Laugwitz's
interpretation and the idea that infinitesimals play a foundational
role in Cauchy.%
\footnote{Felscher's text is examined in more detail in
Subsection~\ref{felscher}.}

In his zeal to do away with Cauchy's infinitesimals, Felscher seeks to
describe them in terms of the modern terminology of {\em germs\/}.
Here two sequences are in the same germ if they agree at infinity,
i.e., for all sufficiently large values of the index~$n$ of the
sequence~$\langle u_n : n\in \N\rangle$, that is, are equal almost
everywhere.%
\footnote{Equality ``almost everywhere'' is a 20th century concept;
Felscher seems to suggest that Cauchy thought of the relation between
his variable quantities and his infinitesimals in a way that would be
later described as equality almost everywhere.}
Felscher's maneuver successfully eliminates the term ``infinitesimal"
from the picture, but has the effect of undermining Felscher's own
thesis, by lending support to the presence of such ``collapsing'' in
Cauchy.

Indeed, the idea of reading germs of sequences into Cauchy is
precisely Laugwitz's thesis in \cite[p.~272]{Lau87}.  Germs of
sequences are also the basis of
Laugwitz's~$\Omega$-calculus~\cite{SL}, a ring constructed using a
Fr\'echet filter.%
\footnote{In the Schmieden--Laugwitz construction, a Fr\'echet filter
is used where an ultrafilter would be used in a hyperreal
construction; see Appendix~\ref{rival}.}

In an editorial footnote to Lakatos's essay, J.~Cleave outlines a
non-Archimedean system developed by Chwistek \cite{Ch}, involving a
quotient by a Fr\'echet filter (similarly to the Schmieden--Laugwitz
construction).  Cleave then concludes that the relation of such
infinitesimals to Cauchy's is ``obvious''.%
\footnote{Cleave writes in his footnote 32*: ``A construction of
non-standard analysis is given in Chwistek \cite{Ch} (1948) which is
derived from a paper published in 1926.  It is basically the reduced
power~$\R^\N/F$ where~$F$ is the Fr\'echet filter on the natural
numbers (the collection of cofinite sets of natural numbers) (see
Frayne, Morel, and Scott \cite{FMS}) [\dots]  This particular
construction is not an elementary extension of~$\R$ but there are
sufficiently powerful transfer properties to enable some non-standard
analysis to be performed.  It may be observed that the elements of
$\R^\N/F$ are equivalence classes of sequences of reals, two sequences
$s_1, s_2 \ldots$ and~$t_1, t_2, \ldots$ being counted equal if for
some~$n$,~$s_m = t_m$ for all~$m \geq n$.  The relation of these
classes to Cauchy's variables is obvious'' \cite[p.~160]{La78}.}
To us it appears that the said relation requires additional argument.

Cauchy's published work contains evidence that he intuitively sensed a
``collapsing" involved in the passage from a null sequence to an
infinitesimal.  Thus, the second chapter of the {\em Cours
d'Analyse\/}~\cite{Ca21} of 1821 contains a series of theorems (eight
of them) whose main purpose appears to be to emphasize the importance
of the {\em asymptotic\/} behavior of the sequence (i.e., as the index
tends to infinity).  This thread is pursued in more detail in
Subsection~\ref{two}.  Furthermore, he refers to his infinitesimals as
``quantities", a term he uses in the context of an ordered number
system, as opposed to the complex numbers which are always
``expressions" but never ``quantities".  In fact, the recent English
translation \cite{BS} of the {\em Cours d'Analyse\/} erroneously
translates one of Cauchy's complex ``expressions" as a ``quantity",
and the reviewer for {\em Zentralblatt\/} dutifully notes this error.%
\footnote{\label{rein1}The reviewer, Reinhard Siegmund-Schultze, was
far more sensitive to Cauchy's infinitesimals on this occasion than in
an earlier instance, occasioned by a review of Felscher's text; see
footnote~\ref{rein2}.}

A second aspect of collapsing is the mental transformation of the
process of ``tending to zero'' into a concept/noun (as a null sequence
is transformed into an infinitesimal), thought of as a reification, or
encapsulation, of a highly compressed process.%
\footnote{In the education literature, such a compressing phenomenon
is studied under the name of procept (process+concept) \cite{GT},
encapsulation~\cite{Du91}, and reification~\cite{Sf91}.}
In the successive theorems in his Chapter~2, Cauchy seeks to collapse
the initial articulation of his infinitesimal~$\alpha$ as a
temporally-deployed {\em process\/}, by deliberately leaving out the
implied index (i.e. label of the terms of the sequence), and by
explicitly specifying and emphasizing a rival index: the exponent in a
power~$\alpha^n$ of the infinitesimal, thought of as an infinitesimal
of a higher and higher order.

We will examine Cauchy's approach further in Section~\ref{two} by
means of a detailed analysis of his text.

\subsection{An 
analysis of {\em Cours d'Analyse\/} and its infinitesimals}
\label{two}

We will refer to the pages in the {\em Cours d'Analyse\/} \cite{Ca21}
itself, rather than the collected works.

Cauchy's Chapter 2, section 1 starts on page 26.  Here Cauchy writes
that a variable quantity becomes {\em infinitely small\/} if, etc.
Here ``infinitely small" is an adjective, and is not used as a
noun-adjective pair.

On page 27, Cauchy employs the noun-adjective combination, by
referring to ``infinitely small quantities" ({\em quantit\'es
infiniment petites\/}, still in the feminine).  He denotes such a
quantity~$\alpha$.  Note that the index in the implied sequence is
suppressed (namely,~$\alpha$ appears without a lower index).%
\footnote{Note that Cauchy uses lower indices to indicate terms in a
sequence in his proof of the intermediate value theorem \cite[Note
III, p.~462]{Ca21}.}

On page 28, he introduces a competing numerical index, namely the
exponent, by forming the infinitesimals
\[
 	\alpha, \alpha^2, \alpha^3, \dots
\]
By the time we reach the bottom of the page (fourth line from the
bottom), he is already employing ``infinitely small" as a {\em noun\/}
in its own right: {\em infiniment petits\/} (in the masculine plural).

On page 29, Theorem 1 asserts that a highest-order infinitesimal will
be smaller than all the others (infinitesimals are consistently
referred to in the masculine).  The theorem has not yet chosen a {\em
letter\/} label for the competing index (i.e. order of infinitesimal).

Still on page 29, Theorem 2 for the first time introduces a label for
the order of the infinitesimal, namely the letter~$n$, as in
\[
\alpha^n.
\]
 
On page 30, Theorem 3 introduces several different letter indices:
\[
 	n, n', n'' \dots
\]
for the orders of his infinitesimals.  The theorem concerns the order
of the sum, again compelling the student to focus on the competing
orders (at the expense of the suppressed index of the ``variable
quantity" itself).

Still on page 30, Theorem 4 introduces the terminology of {\em
polynomials\/} in~$\alpha$, and describes their orders~$n, n', n''
\ldots$ for the first time as a ``sequence".  We now have two
``sequences": (the variable quantity)~$\alpha$ itself, whose index is
implicit (was never labeled), and the sequence of orders, which are
both emphasized and elaborately labeled using ``primes''~$'$ and
double primes~$''$.

In all these theorems, it is the {\em asymptotic\/} behavior of null
sequences that is constantly emphasized, which suggests that Cauchy
might have found it perfectly natural to identify/collapse sequences
that agree almost everywhere.  Terminology {\em finit par \^etre\/},
{\em finit par devenir\/} (suggestive of such collapse) is employed
repeatedly.

On pages 31 and 32, three additional theorems and one corollary are
stated, for a total of eight results on the asymptotic behavior of
such null sequences.

\renewcommand{\arraystretch}{1.3}
\begin{table}
\[
\begin{tabular}[t]
{ | p{1.2in} || p{1.5in} | p{1.4in} | p{.5in} | p{.75in} | p{.5in} |}
\hline & independent variable increment ($\Delta x$) & dependent
variable increment ($\Delta y$) \\ \hline\hline Cauchy's first
definition & infinitesimal & variable tending to zero \\ \hline
Cauchy's second definition & infinitesimal & infinitesimal \\ \hline
\end{tabular}
\]
\caption{\textsf{Cauchy's first two definitions of continuity in 1821
are of the form {\em ``if~$\Delta x$ is~$\ldots$, then~$\Delta y$ is
$\ldots$''\/}.  Note the prevalence of the term ``infinitesimal''.}}
\label{continuity}
\end{table}
\renewcommand{\arraystretch}{1}

By the time Cauchy reaches Section 2 of Chapter 2 on page 34
(concerning continuity of functions), he has already encapsulated the
{\em process\/} implied in the notion of ``variable quantity", into
the concept/masculine noun {\em infiniment petit\/}.  When he evokes
an infinitely small~$x$-increment~$\alpha$, only a stubborn
Weierstrassian will refuse to interpret his~$\alpha$ as a
concept/noun.  The {\em second\/} definition (out of the three
definitions of continuity given here) is the one Cauchy italicizes,
implying it is the main one.  Here both the~$x$-increment and
the~$y$-increment are described as infinitely small increments (see
Table~\ref{continuity}).

In this context, the verb ``become" is being used in two different
senses: 
\begin{enumerate}
\item[(a)] the terms in the sequence {\em become\/} smaller than any
number;
\item[(b)] the encapsulating sense of a process being compressed into
(and thus {\em becoming\/}) a concept/noun,
\end{enumerate}
as analyzed by Sad, Teixeira, and Baldino \cite{STB}, who employed the
terminology of a {\em transformation of essence\/}.  It is interesting
to note that Bolzano fought against the sense (a), by suppressing the
parametrisation altogether, and viewing the null sequence as a {\em
set\/} (perhaps he was influenced by Zeno paradoxes), but not against
the sense (b).

Note Cauchy's emphasis on the {\em noun\/} aspect of his
infinitesimals:
\begin{quote}
Lorsque les valeurs num\'eriques%
\footnote{\label{num}The meaning of the expression {\em valeur
num\'erique\/} is subject to debate; see next section and
footnote~\ref{sch}.}
successives d'une m\^eme variable d\'ecroissent ind\'efiniment, de
mani\`ere \`a s'abaisser au-dessous de tout nombre donn\'e, cette
variable devient ce qu'on nomme {\em un infiniment petit\/} ou une
quantit\'e infiniment petite.  Une variable de cette esp\`ece a z\'ero
pour limite \cite[p.~4]{Ca21}.
\end{quote}
The use of the noun, {\em un infiniment petit\/}, makes it difficult
to interpret the ``becoming'' in the sense (a) above; rather, the
definition requires sense (b) to be grammatically coherent.  Here the
variable [quantity] {\em becomes\/} a masculine noun: {\em un
infiniment petit\/}.  Cauchy is very precise here: it is the {\em
limit\/} of the variable that's zero.  The variable itself {\em
becomes\/} an {\em infiniment petit\/}.  Cauchy wrote neither that a
variable {\em is\/} an infinitesimal, nor that the limit of the
infinitesimal is zero, but rather that the limit of the {\em
variable\/} is zero, cf.~\cite[p.~301-302]{STB}.

Once the use of {\em infiniment petit\/} as a noun is established,
Cauchy freely uses the term interchangeably as a noun or as an
adjective.

Cauchy's {\em Cours d'Analyse\/} presented only a theory of
infinitesimals of polynomial rate of growth as compared to a
given~$\alpha$.  His theory of infinitesimals of more general order,
and its influence on later authors such as E.~Borel, was already
discussed in Section~\ref{one};%
\footnote{See also the main text following footnote~\ref{ehr}.}
see Section~\ref{time} for a broader
historical perspective.  In Subsection~\ref{micro1}, we will discuss a
specific application Cauchy makes of his infinitesimals, in the
context of the sum theorem.

\subsection{Microcontinuity and Cauchy's sum theorem}
\label{micro1}

How is Cauchy able to define concepts such as uniform continuity and
uniform convergence in terms of a {\em single\/} variable, unlike
standard definitions thereof which call for a {\em pair\/} of
variables?

Let~$x$ be in the domain of a function~$f$, and consider the following
condition, which we will call {\em microcontinuity\/} at~$x$:
\begin{quote}
``if~$x'$ is in the domain of~$f$ and~$x'$ is infinitely close to~$x$,
then~$f(x')$ is infinitely close to~$f(x)$".
\end{quote}
Then ordinary continuity of~$f$ is equivalent to~$f$ being
microcontinuous on the Archimedean continuum (A-continuum for short),
i.e., at every point~$x$ of its domain in the A-continuum.  Meanwhile,
uniform continuity of~$f$ is equivalent to~$f$ being microcontinuous
on the Bernoulian continuum (B-continuum for short), i.e., at every
point~$x$ of its domain in the B-continuum (the relation of the two
continua is discussed in more detail in Appendix~\ref{rival}).

Consider, for example, the function~$\sin(1/x)$ defined for
positive~$x$.  The function fails to be uniformly continuous because
microcontinuity fails at a positive infinitesimal~$x$ (due to ever
more rapid oscillation).  The function~$x^2$ fails to be uniformly
continuous because of the failure of microcontinuity at an infinite
member of the B-continuum.

A similar distinction exists between pointwise convergence and uniform
convergence.%
\footnote{See e.g.,~Goldblatt \cite[Theorem 7.12.2, p.~87]{Go}.}
Which condition did Cauchy have in mind in 1821?  Abel interpreted it
as convergence on the A-continuum, and presented ``exceptions'' (what
we would call today counterexamples) in 1826.  Additional such
exceptions were published by Seidel and Stokes in the 1840s.

Cauchy's contemporary E. G. Bj\"orling anticipated the B-continuum
definition in a remarkable passage dating from 1852:
\begin{quote}
when one comes to show that in a series, of which the terms are
functions of a quantity~$x$, converges [\ldots] for each given value
of~$x$ up to a certain limit~$X$, it is not necessary to believe that
the series continues necessarily to converge [\dots] for values of~$x$
indefinitely close to that limit (Bj\"orling \cite[p.~455]{Bj52a}, as
cited by Grattan-Guinness \cite[p.~232]{Grat86}).
\end{quote}
Here the ``limit''~$X$ is clearly a member of the usual Archimedean
continuum; but what is the nature of the~$x$ mentioned at the end of
Bj\"orling's phrase?  Bjorling's reference to convergence at~$x$
implies that~$x$ is understood to be an individual/atomic element of
the domain; yet speaking of~$x$ as being ``indefinitely close'' to~$X$
is meaningless in the context of an A-continuum.  One faces a stark
choice of either toeing the triumvirate line on the ``true'' continuum
and therefore finding Bj\"orling's arguments ``absurd'', as Pringsheim
did \cite[p.~345]{Pr};%
\footnote{As cited by Grattan-Guinness \cite[p.~233]{Grat86}.}
or, removing the blinders so as to envision a mid-19th century
anticipation of an infinitesimal-enriched continuum.

Cauchy clarified/modified his position on the sum theorem in 1853.%
\footnote{Grattan-Guinness \cite{Grat86} and Br\aa ting \cite{Br07}
argue that Cauchy was influenced by Bj\"orling.}
In his text \cite{Ca53}, he specified a stronger condition of
convergence on the B-continuum, including at~$x=\frac{1}{n}$
(explicitly mentioned by Cauchy).  The stronger condition bars Abel's
counterexample.

To give a more detailed explanation of Cauchy's 1853 text, note that
Cauchy's approach is based on two assumptions which can be stated in
modern terminology as follows:
\begin{enumerate}
\item
when you have a closed expression for a function, then its values at
``variable quantities'' (such as $x=\frac{1}{n}$) are calculated by
using the same closed expression as at real values;
\item
to evaluate a function at a variable quantity generated by a sequence,
one evaluates term-by-term.
\end{enumerate}
Now in 1853 Cauchy analyzed Abel's counterexample (without mentioning
Abel's name) by first writing down the closed form of the sum of the
{\em remainder\/} of the series.  This is given by a certain integral.
He proceeded to evaluate it at an infinitesimal $x=\frac{1}{n}$ using
assumption~(1).  Concretely, he substituted the values $\frac{1}{n}$
into the closed expression (integral) using assumption~(2).  The
sequence he got was {\em not\/} a null sequence.  He concluded that
the remainder term at this particular infinitesimal is not an
infinitesimal.  Hence, he concluded, the ``always'' part of his
hypothesis is not satisfied by Abel's example.

Thus, Cauchy's clarification/modification from 1853 amounts to
requiring convergence on an incipient form of a B-continuum.  Some
historians acknowledge his clarification/modification, and interpret
it as the addition of the condition of uniform convergence, while
adhering to an A-continuum framework.  The recent text \cite{KK11b}
argued that such an interpretation is problematic.  Namely, Cauchy
states the condition in terms of a {\em single\/} variable, whereas
the traditional definition of uniform continuity or convergence in the
context of an A-continuum necessarily requires {\em a pair\/} of
variables.  Cauchy specifically evokes~$x=\frac{1}{n}$ where Abel's
``exception/counterexample" fails to converge.  The matter is
discussed in detail by K.~Br\aa ting \cite{Br07}.

\subsection{Five common misconceptions in the Cauchy literature}
\label{25}

Both Laugwitz and Hourya Benis Sinaceur \cite{Si} have exposed a number
of misconceptions in the literature concerning Cauchy's foundational
work.  Five of the most common ones are reproduced below in italics.

\medskip\noindent 1. {\em Bolzano, Cauchy, and Weierstrass were all
  gardeners who contributed to the ripening of the fruit of the notion
  of limit.\/}

\medskip
Here the implicit assumption is that the Weierstrassian epsilontic
notion of ``limit" in the context of an Archimedean continuum is the
centerpiece of any possible edifice of analysis.  Such an assumption
is questionable on two counts.  First, as Felix Klein pointed out in
1908, there are two parallel threads in the development of analysis,
one based on an Archimedean continuum, and the other exploiting an
infinitesimal-enriched continuum.%
\footnote{See material in the main text around footnote~\ref{klein2}.}
One risks pre-judging the outcome of any analysis of Cauchy by
postulating that he is working in the Archimedean thread. The second
implicit assumption is that the Weierstrassian notion of limit is
central in Cauchy.  Thus, Boyer \cite[p.~277]{Boy} postulates that
Cauchy is working with a notion of limit similar to the Weierstrassian
one.  This requires further argument, and at least at first glance is
incorrect: Cauchy emphasizes infinitesimals as a foundational notion,
but he never emphasizes limits as a foundational notion. Thus, in
Cauchy's definition of continuity the word ``limit" does occur, but
only in the sense of the ``endpoint" of the interval of definition of
the function, rather than the behavior of its values.%
\footnote{The misconception may be found, for instance, in J. Gray,
see main text around footnote~\ref{gray}.}

\medskip\noindent
2. {\em Cauchy, along with other mathematicians, abandoned
    infinitesimals in favor of other more rigorous notions.\/} 

\medskip
During the period 1814-1820, Cauchy appears to have been ambivalent
about infinitesimals.  Starting in about 1821, he uses them with
increasing frequency both in his textbooks and his research
publications, and insists on the centrality of infinitesimals as a
foundational notion; see Section~\ref{time} for a related chronology.

\medskip\noindent
3. {\em Cauchy was forced to teach infinitesimals at the Ecole.\/}

\medskip
The intended implication appears to be that Cauchy only used
infinitesimals because of the pressure from the Ecole
administration. During the period 1814-1820 there were some tensions
with the administration over the delayed appearance of infinitesimals
in the syllabus, see \cite{Gil}.  At any rate, Cauchy continued using
infinitesimals thoughout his career and long after completing his
teaching stint at the Ecole in 1830.  Thus he reproduces his 1821
definition of continuity (in terms of infinitesimals) as late as 1853,
in his text on the sum theorem \cite{Ca53}, see
Subsection~\ref{micro1}.

\medskip\noindent
4. {\em Cauchy based his infinitesimals on the notion of limit.\/}

\medskip
This is an ambiguous claim, and essentially a play on words on the
term ``limit".  The modern audience understands ``limit" as a
Weierstrassian epsilontic notion.  If this is what is claimed, then
the claim is false.  As far as the kinetic notion of limit that Cauchy
does mention in discussing a variable quantity approaching a limit, it
is conspicuously absent in Cauchy's discussion of infinitesimals.
Thus, rather than infinitesimals being based on the notion of limit,
it is the notion of a variable quantity that's primitive, and both
infinitesimals and limits are defined in terms of it (see
Subsection~\ref{gilain} for a more detailed discussion).  Note again
that the term ``limit" does appear in Cauchy's definition of
continuity, but in an entirely different sense, namely endpoint of the
interval where the function is defined.%
\footnote{\label{ult}To elaborate, we will note that the purported
opposition between infinitesimals and limits (a frequently found claim
in the literature) is a conceptual error in its own right.  The
opposition is not between limits and infinitesimals.  Limits are
present in both approaches.  In the infinitesimal approach, limits can
be defined in terms of the standard part function, and the latter in
terms of limits.  The true opposition is between infinitesimals, on
the one hand, and epsilon-delta with its quantifier complications, on
the other.  The preoccupation with the word ``limit'' can be exposed
as a purely linguistic one by the following thought experiment, based
on Pourciau's lucid analysis, see footnote~\ref{pourciau}.  Imagine
that Newton had developed an ``ult'' notation for his ultimate ratios,
e.g.,~$\text{ult}_{h\to 0}\frac{f(x+h)-f(x)}{h}$ for the derivative,
and~$\text{ult}_{n\to \infty}\sum_{i=1}^n f(x_i)\Delta x$ for the
integral.  Post-Weierstrassian historians would then perhaps be more
moderate in their enthusiasm for the vague musings about limits
(incorrectly) attributed to d'Alembert, see
Subsection~\ref{chapelle}.}

\medskip\noindent
5. {\em Cauchy introduced rigor into calculus that anticipates the
    rigor of Weierstrass.\/}

\medskip
While Cauchy certainly emphasizes rigor, postulating a continuity
between Cauchy's rigor and Weierstrassian rigor is a methodological
error.%
\footnote{See Klein's comments discussed in the main text around
footnote~\ref{klein2}, concerning the two parallel strands for the
development of analysis.}
To Cauchy, rigor meant abandoning the principle of the ``generality of
algebra" as practiced by Euler, Lagrange, and others, and its
replacement by geometry--and by infinitesimals.

\section{Comments by modern authors}
\label{modern}

In this section, we analyze some comments by modern authors related to
Cauchy's foundational contribution.

\subsection{Gilain's limit}
\label{gilain}

C.~Gilain affirms the following:
\begin{quote}
On sait que Cauchy d\'efinissait le concept d'infiniment petit \`a
l'aide du concept de limite, qui avait le premier r\^ole%
\footnote{``We know that Cauchy defined the concept of infinitely
small by means of the concept of limit, which played the primary
role.''}  (voir Analyse alg\'ebrique, p.~19 \dots)''
\cite[footnote~67]{Gil}.
\end{quote}
Here Gilain is referring to Cauchy's Collected Works, S\'erie 2,
Tome~3, p.~19, corresponding to \cite[p.~4]{Ca21}.  Both of Gilain's
claims are erroneous, as we now show.  Cauchy starts by discussing
{\em variable quantities\/} as a primary notion, in the following
terms: 
\begin{quote}
On nomme quantit\'e {\em variable\/} celle que l'on consid\`ere comme
devant re\c cevoir successivement plusieurs valeurs diff\'erentes les
unes des autres.
\end{quote}
Next, Cauchy exploits his primary notion to evoke his kinetic concept
of limit%
\footnote{Such a concept is similar to Newton's, see
footnote~\ref{pourciau}.}
in the following terms: 
\begin{quote}
Lorsque les valeurs successivement attribu\'ees \`a une m\^eme
variable s'approchent ind\'efiniment d'une valeur fixe, de mani\`ere
\`a finir par en diff\'erer par aussi peu qu l'on voudra, cette
derni\`ere est appel\'ee la {\em limite\/} de toutes les autres.
\end{quote}
Finally, Cauchy proceeds to define infinitesimals in the following
terms: 
\begin{quote}
Lorsque les valeurs num\'eriques successives d'une m\^eme variable
d\'ecroissent ind\'efiniment, de mani\`ere \`a s'abaisser au-dessous
de tout nombre donn\'e, cette variable devient ce qu'on nomme {\em un
infiniment petit\/} ou une quantit\'e infiniment petite.  Une variable
de cette esp\`ece a z\'ero pour limite.%
\footnote{This definition was already analyzed in
Subsection~\ref{two}, see main text around footnote~\ref{num}.}
\end{quote}
Thus, Cauchy defined both infinitesimals and limits in terms of
variable quantities.  Neither is the limit concept primary, nor are
infinitesimals defined in terms of limits, contrary to Gilain's
claims.%
\footnote{For a similar misconception in Smithies, see
footnote~\ref{smi}.}

\subsection{Felscher's {\em Bestiarium infinitesimale\/}}
\label{felscher}

We have argued for an interpretation of Cauchy's foundational stance
that endeavors to take Cauchy's infinitesimals at their face value.
Such an interpretation has not been without its detractors.  A decade
ago, W.~Felscher \cite{Fe} set out to investigate Cauchy's continuity,
in an 18-page text, marred by an odd focus on d'Alembert.%
\footnote{Felix Klein \cite[p.~103]{Kl} discusses the error in
d'Alembert's proof of the fundamental theorem of algebra, first
noticed by Gauss.}
To be sure, it is both legitimate and necessary to examine Cauchy's
predecessors, including d'Alembert, if one wishes to understand Cauchy
himself.  Indeed, a debate of long standing (over a century long, in
fact) had opposed two rival methodologies in the study of the
foundations of the new science of Newton and Leibniz:
\begin{enumerate}
\item[(A)] a methodology eschewing infinitesimals; and
\item[(B)] a methodology favoring them.%
\footnote{See Felix Klein's comments discussed in the main text around
footnote~\ref{klein2}.}
\end{enumerate}
It is legitimate to ask which of the two methodologies is the one that
underpins Cauchy's oeuvre.  However, Felscher's conceptual framework
is flawed in a fundamental way.  The outcome of his investigation is
predetermined from the outset by the following two factors:
\begin{enumerate}
\item
Felscher's exclusive focus on d'Alembert,%
\footnote{Felscher mentions Euler and the Bernoullis in his section
entitled {\em D'Alembert's program\/}, but says not a word about them
in his article.}
one of the radical adherents of the A-methodology, and 
\item
Felscher's {\em postulating\/} a methodological continuity between the
work of d'Alembert and Cauchy (see also Subsection~\ref{chapelle}).
\end{enumerate}
Displaying a masterly command of scholarly Latin, Felscher offers the
reader a glimpse of the {\em bestiarium infinitesimale\/} in section~6
of his essay, starting on page~856.  The punchline comes in the middle
of page~857, where Felscher points out that Cauchy refers specifically
to {\em numerical\/} values of his variables, the latter being
described by Cauchy as {\em becoming\/} infinitesimals.

The adjective {\em numerical\/} is linked etymologically to the noun
{\em number\/}.  Cauchy's numbers (unlike his {\em quantities\/}) are
certainly appreciable (i.e., neither infinitesimal nor infinite, nor
even negative), as can be seen by reading the first page {\em
Pr\'eliminaires\/} of his 1821 {\em Cours D'Analyse\/}.  Felscher
concludes that the variables assume only appreciable values, but not
non-Archimedean ones.%
\footnote{\label{sch}Note that Schubring \cite[p.~446]{Sch}, in
footnote 14, explains Cauchy's term {\em numerical value\/} as what we
would call today the {\em absolute value\/}.  Fisher \cite[p.~262]{Fi}
interprets Cauchy's definition accordingly, so as to allow room for
infinitesimal values of Cauchy's variables.  See also
footnote~\ref{Luz2} on Cleave.}

Felscher's etymological insight offers a refutation of the Luzin
hypothesis,%
\footnote{Luzin himself, in fact, similarly rejected non-Archimedean
time, as discussed by Medvedev \cite{Me93}.}
to the effect that Cauchy variables may pass through non-Archimedean
values on their way to zero.  Has Felscher shown that Cauchy's {\em
bestiarium infinitesimale\/} is in fact uninhabited?

Hardly so.  While chasing out the infinitesimal mouse of Luzin's
hypothesis%
\footnote{\label{Luz2}Luzin was probably not the first and surely not
the last to formulate a hypothesis to the effect that Cauchy's
variable quantities pass through infinitesimal values on their way to
zero.  Lakatos \cite[p.~153]{La78} speculates that Cauchy's variables
``ran through Weierstrassian real numbers {\em and\/}
infinitesimals'', while J.~Cleave (who edited Lakatos's essay for
publication in the {\em Mathematical Intelligencer\/}) in footnote 18*
in \cite[p.~159]{La78} disagrees, limiting Cauchy variables to
sequences of Weierstrassian reals (Cleave quotes the relevant passage
on {\em numerical values\/} but does not analyze it here).  Cleave
alludes to the etymological point in \cite[p.~268]{Cl79}, where he
disagrees with Fisher on this point.}
for the Cauchy variable, Felscher missed the elephant of the
possibility of the variable {\em becoming\/} an infinitesimal through
a process of reification (see Subsection~\ref{second}).

On page~846, Felscher quotes an agitated passage from d'Alembert's
1754 article.  D'Alembert attacks the {\em obscurity\/}, and even the
{\em falsehood\/} of a definition of infinitesimals attributed to
unnamed geometers, and sums up his thesis by accusing such geometers
of {\em charlatanerie\/}, a term ably translated as {\em quackery\/}
by Felscher, who sums up as follows:
\begin{quote}
Reading these words today we may get the impression that they were
written at the time of Weierstrass or Cantor,%
\footnote{Cantor was indeed a worthy heir to d'Alembert's
anti-infinitesimal vitriol.  Cantor dubbed infinitesimals the {\em
cholera bacillus\/} of mathematics, see J.~Dauben \cite[p.~353]{Da95},
\cite[p.~124]{Da96}, H.~Meschkowski \cite[p.~505]{Me}.  This was
perhaps the most vitriolic opposition to the B-continuum (see
Appendix~\ref{rival}) prior to Errett Bishop's {\em debasement of
meaning\/}, a term he applied to classical mathematics in general in
1973 \cite{Bi85}, and to infinitesimal calculus \`a la Robinson in
particular, in 1975 \cite{Bi75}, see \cite{KK11a} and \cite{KK11d} for
details.}
or even by a contemporary mathematician.%
\footnote{It is worth pondering which contemporary mathematician
(known for anti-infinitesimal vitriol) Felscher may have had in mind
here, given his interest in intuitionistic logics \cite{Fe85, Fe86},
see also \cite{KK11d}.}
\end{quote}
It is sobering to realize that, forty years after A.~Robinson, a
logician named Walter Felscher still conceived of the history of
analysis in terms of a triumphant march out of the dark age of the
infinitesimal, and toward the yawning heights of Weierstrassian
epsilontics.%
%
%

D'Alembert's verbal excesses merely put in relief the fact that no
such rhetoric is to be found anywhere in either Cauchy or Bolzano.
Felscher presents a convincing case that d'Alembert was opposed to
infinitesimals.  Felscher's title {\em Bolzano, Cauchy, epsilon,
delta\/}%
\footnote{Apparently, a kind of a mantra: {\em Bolzano, Cauchy,
epsilon, delta; Bolzano, Cauchy, epsilon, delta; \dots\/} which,
repeated sufficiently many times, would lead one to accept Felscher's
reduction of Cauchy's continuum to an A-continuum (see
Appendix~\ref{rival}).}
could therefore have pertinently been replaced by {\em D'Alembert,
Weierstrass, epsilon, delta\/}, as the case for Bolzano's opposition
to infinitesimals can similarly be challenged.%
\footnote{In discussing Bolzano's attitude toward infinitesimals, we
have to distinguish between the early Bolzano and the late Bolzano.
The early Bolzano defines the ``Infinitely small" as ``variable
quantities" in the following terms: A quantity is infinitely small if
it becomes less than any given quantity (here Bolzano does not speak
of ``values" [but naturally he thought of them]).  The late Bolzano
definies infinitely small (and infinitely large) {\em numbers\/}; one
of them is~$1/(1+1+1+\ldots)$ (infinitely many terms).  We are
grateful to D.~Spalt for this historical clarification.}
Indeed, as Lakatos points out,
\begin{quote}
[Bolzano] was possibly the only one to see the problems related to the
difference between the two continuums: the rich Leibnizian continuuum
and, as he called it, its `measurable' subset--the set of
Weierstrassian real numbers.  Bolzano makes it very clear that the
field of `measurable numbers'%
\footnote{\label{measure}{\em Measurable number\/} is Bolzano's term
for appreciable number (no relation to Lebesgue-measurability).  Here
Bolzano foreshadows Bj\"orling's dichotomy (see \cite{Br07}), which can
be analyzed in terms of A- and B-continua (see Appendix~\ref{rival}).}
constitutes only an Archimedean subset of a continuum enriched by
non-measurable - infinitely small or infinitely large - quantities
\cite[p.~154]{La78}.
\end{quote}

D.~Kurepa \cite[p.~664]{Ku} provides some details on Bolzano's use of
infinitesimals.  As far as epsilon-delta techniques are concerned, a
case for Bolzano's anticipation thereof is far more convincing than
for Cauchy.%
\footnote{Hourya Benis Sinaceur wrote: ``les oeuvres de Bolzano et de
Cauchy repr\'esentent deux courants distincts, h\'et\'erogen\`es \`a
l'origine, et qui ne se rencontrent pas avant les ann\'ees 1870, avec
les travaux de Weierstrass, Cantor, etc.  Seule l'illusion
r\'etrospective et une connaissance indirecte des textes permettent de
les amalgamer dans un m\^eme cours ininterrompu suppos\'e traverser
tout le XIX$^e$ si\`ecle'' \cite[p.~111]{Si}, and cites a letter of
H. A. Schwarz to Cantor \cite{Sc} to the effect that ``la m\'ethode de
d\'emonstration [\ldots] de Weierstrass est un d\'eveloppement des
principes de Bolzano'' \cite[p.~112]{Si}.}

Felscher's intriguing parenthetical remark indicates that he was more
sensitive to Cauchy's language than numerous Cauchy historians:
\begin{quote}
it is left open whether a {\em quantit\'e variable\/}, with an
assignment converging to zero, actually {\em is\/} or only {\em
becomes\/} a {\em quantit\'e infiniment petite\/} \cite[p.~850]{Fe}.
\end{quote}

On page 851, Felscher presents an analysis of Cauchy's use of an
infinitesimal quantity, denoted~$i$, in differentiating an exponential
function.  Here Felscher's additional parenthetical remark, to the
effect that ``notational confusion arises from denoting both the
variable~$i$ and its values by the same letter'', is an unjustified
criticism of Cauchy.  The criticism underscores Felscher's
insensitivity toward the dynamic aspect of Cauchy's infinitesimal~$i$,
when individual values are irrelevant in the context of the dynamism
of the encapsulation taking place whenever Cauchy evokes an
infinitesimal.  Otherwise Felscher's analysis is unexceptionable, save
for a {\em non-sequitur\/} of a conclusion:
\begin{quote}
``No `infinitesimal' non-Archimedean numbers are ever used by Cauchy
for his {\em quantit\'es infiniment petites\/}.''
\end{quote}
In reality, Cauchy's discussion of the derivative of the exponential
function admits a number of possible interpretations.

On page 852, Felscher analyzes Cauchy's infinitesimals in modern
terms:
\begin{quote}
Using today's terminology, one would describe Cauchy's forms to be
filled by assignments as functions, but in order to distinguish them
from the actual functions subsequently considered by Cauchy, one might
call them {\em functional germs\/}. [Emphasis in the
original--authors]
\end{quote}
Felscher mentions Cauchy's use of functional germs again on page 855.
In his zeal to rename Cauchy's infinitesimals by employing a modern
notion, so as bashfully to avoid the distasteful {\em infi\/} term,
Felscher comes close to endorsing Laugwitz's ``Cauchy numbers'',
similarly defined in terms of germs.

Thus, Laugwitz wrote: 
\begin{quote}
Every real function~$f(u)$ defined on an interval~$0 < u < p$
represents a Cauchy number.  Two such functions~$f(u)$ and~$g(u)$
represent the same Cauchy number if and only if there is an interval
$0 < u < q$ in which~$f(u) = g(u)$ \cite[p.~659]{Lau97}.
\end{quote}
Note that Laugwitz essentially defines his ``Cauchy numbers'' by
exploiting the concept of the germ of a function; Laugwitz explicitly
mentions {\em function germs\/} in \cite[p.~272]{Lau87}.  Felscher
\cite[p.~858]{Fe} leaves very little doubt as to how he felt with
regard to Laugwitz infinitesimals:
\begin{quote}
In this connection one must also mention certain articles and books by
D.~Laugwitz, in which [\dots] he develops his own `mathematics of the
infinitesimal' and uses it to interpret skillfully various aspects of
the mathematics of the period from Euler to Cauchy.  And so we have
glanced at the {\em bestiarium infinitesimale\/}.%
\footnote{The reviewer for MathSciNet concurred, see
Section~\ref{belluot}.}
\end{quote}

On page 853, Felscher {\em omits\/} a crucial first phrase used by
Cauchy in formulating his first definition of continuity.  Namely, he
omits Cauchy's phrase {\em stating that~$\alpha$ is an
infinitesimal\/}, corresponding to the upper-left entry in
Table~\ref{continuity}.%
\footnote{\label{crucial}This crucial detail leads Felscher to a
further error of a conceptual nature, discussed below.}

On pages 854-855, Felscher makes the following statement concerning
Cauchy's notion of {\em continuity\/}:
\begin{quote}
Both Bolzano and Cauchy gave definitions of continuity which express
today's [\ldots] continuity.  Both made their definitions precise and
used them in today's sense; both employed them by comparing numbers
and their distances with the help of inequalities in order to prove
important theorems in analysis.  However, Cauchy defined and used the
notion of limit, whereas Bolzano did not.
\end{quote}
Felscher's assertion concerning continuity and inequalities is
misleading.  Generally speaking, Cauchy's limit concept is a kinetic
one (see Hourya Benis Sinaceur~\cite{Si}), and is a derived notion,
depending for its definition on the primary notion of a {\em variable
quantity\/}.  On occasion, Cauchy worked with real inequalities in
proofs (see, e.g., Grabiner \cite{Grab}), but he never gave such a
definition of {\em continuity\/}.

On page 855, line 9, Felscher alleges that Cauchy's first definition
of continuity is similar to Bolzano's, with the implication that the
terminology of ``infinitesimal'' is not employed by Cauchy.  Now the
form in which Cauchy's definition was quoted by {\em Felscher\/} two
pages earlier (see our comment above concerning Felscher's page 853)
did not employ infinitesimals.  But the form in which it appears in
{\em Cauchy\/} did employ infinitesimals.%
\footnote{See Table~\ref{continuity} (for a summary of Cauchy's
definitions) and footnote~\ref{crucial}.  Schubring \cite[p.~465]{Sch}
writes that J.~L\"utzen's is the best analysis of continuity in
Cauchy.  Meanwhile, L\"utzen \cite[p.~166]{Lut03} states: ``Cauchy
[\ldots] gives two definitions, first one without infinitesimals, and
then one using infinitesimals.''  The second claim is correct, but not
the first.}

The most remarkable aspect of Felscher's, unfortunately seriously
flawed, essay is how close he comes to sensing the cognitive view of
compression/encapsulation outlined in Subsection~\ref{second}:
\begin{quote}
We speak of a function or a variable approaching some value {\em
ind\'efiniment\/} (indefinitely); we imagine a limiting process.
[\ldots] Thus far,~$\epsilon$ and~$\delta$ (and in case of sequences
also~$n$ and~$N$) appear as handles affixed to the stages of those
infinite processes.  It seems that if appropriately handled in our
mental exercises, they enable us to use finitely many arguments to
prove statements that, in the end, speak about all the stages of the
infinite process \cite[p.~858]{Fe}.
\end{quote}

\subsection{D'Alembert or de La Chapelle?}
\label{chapelle}

Felscher's text analyzed in Section~\ref{felscher} describes
d'Alembert as ``one of the mathematicians representing the heroic age
of calculus'' \cite[p.~845]{Fe}.  Felscher buttresses his claim by a
lengthy quotation concerning the definition of the limit concept, from
the article {\em Limite\/} from the {\em Encyclop\'edie ou
Dictionnaire Raisonn\'e des Sciences, des Arts et des M\'etiers\/}
(volume 9 from 1765):
\begin{quote}
On dit qu'une grandeur est la limite d'une autre grandeur, quand la
seconde peut approcher de la premi\`ere plus pr\`es que d'une grandeur
donn\'ee, si petite qu'on la puisse supposer, sans pourtant que la
grandeur qui approche, puisse jamais surpasser la grandeur dont elle
approche; ensorte que la diff\'erence d'une pareille quantit\'e \`a sa
limite est absolument inassignable (Encyclop\'edie, volume 9, page
542).
\end{quote}
One recognizes here a kinetic definition of limit already exploited by
I.~Newton.%
\footnote{See footnote~\ref{pourciau} on Pourciau's analysis.}
Whatever the merits of attributing visionary status to this quote,
what Felscher overlooked is the fact that the article {\em Limite\/}
was written by two authors.  In reality, the above passage defining
the concept of ``limit" (as well as the two propositions on limits)
did not originate with d'Alembert, but rather with the encyclopedist
Jean-Baptiste de La Chapelle.  De la Chapelle was recruited by
d'Alembert to write 270 articles for the {\em Encyclop\'edie\/}.

The section of the article contaning these items is signed (E) (at
bottom of first column of page 542), known to be de La Chapelle's
``signature'' in the {\em Encyclopedie\/}.  Felscher had already
committed a similar error of attributing de la Chapelle's work to
d'Alembert, in his 1979 work \cite{Fe79}.%
\footnote{We are grateful to D.~Spalt for this historical
clarification.}
Note that Robinson \cite[p.~267]{Ro66} similarly misattributes this
passage to d'Alembert.

\subsection{Brillou\"et-Belluot: epsilon-delta, period}
\label{belluot}

An instructive case study in a {\em bestiarium\/}-consignment attitude
toward infinitesimals is the review of Felscher's text for Math
Reviews, by one Nicole Brillou\"et-Belluot \cite{Bri}.  The review
contains not an inkling of the fact that the text in question is a
broadside attack on scholars attempting to analyze Cauchy's
infinitesimals seriously.  The reviewer mentions the ``epsilon-delta
technique", and notes that Cauchy made his ``definitions [of
continuity] precise and used them in today's sense", but fails to
mention that the definitions in question are {\em infinitesimal\/}
ones, a fact not denied by Felscher (at least in the case of one of
the definitions).  Her review mentions d'Alembert, who does not appear
in Felscher's summary, indicating that she had read the body of
Felscher's text itself (rather than merely Felscher's summary).

Brillou\"et-Belluot notes that Felscher reports on how ``limit was
explained and used by d'Alembert and Cauchy".  She reports neither on
Felscher's extensive, and vitriolic, quotes from d'Alembert (including
the dramatic phrase ``the metaphysics and the infinitely small
quantities, whether larger or smaller than one another, are totally
useless in the differential calculus"), nor d'Alembert's colorful
epithets like {\em charlatanerie\/}/quackery.  The reviewer identified
with Felscher's conclusions to such an extent that she chose to spare
the Math Reviews reader the burden of infinitesimal quackery, judging
that Felscher carried that burden once and for all, for the rest of
us.%
\footnote{\label{rein2}Similar remarks apply to the review of
Felscher's text by Reinhard Siegmund-Schultze~\cite{Sie} for
Zentralblatt Math.  Siegmund-Schultze was far less myopic in a later
instance occasioned by a review of an English translation of the {\em
Cours d'analyse\/}, see footnote~\ref{rein1}.}

\subsection{Bos on preliminary explanation}

Leibniz historian H.~Bos acknowledged that Robinson's hyperreals
provide~a
\begin{quote}
preliminary explanation of why the calculus could develop on the
insecure foundation of the acceptance of infinitely small and
infinitely large quantities \cite[p.~13]{Bos}.
\end{quote}

\subsection{Medvedev's delicate question}

Reviewer Cooke \cite{Co} notes that Medvedev \cite{Me87} ``devotes
considerable space to refuting the point of view of modern analysis as
a thing waiting to be born, which earlier mathematicians were striving
unsuccessfully to find and for lack of which they were merely groping
in the dark.  In particular he argues that the approaches of Cauchy
and Weierstrass were so different from each other that it would be
more accurate to speak of two systems of analysis.''

F.~Medvedev further points out that nonstandard analysis
\begin{quote}
makes it possible to answer a delicate question bound up
with earlier approaches to the history of classical analysis.  If
infinitely small and infinitely large magnitudes are [to be] regarded
as inconsistent notions, how could they [have] serve[d] as a basis for
the construction of so [magnificent] an edifice of one of the most
important mathematical disciplines?  \cite{Me87, Me98} 
\end{quote}
A powerful question, indeed.  How do historians answer Medvedev's
question?

\subsection{Grabiner's ``deep insight''}

Not all scholars are satisfied with the {\em
amazing-intuition-and-deep-insight\/} answer offered by J.~Grabiner
who writes:
\begin{quote}
[M]athematicians like Euler and Laplace had a deep insight into the
basic properties of the concepts of the calculus, and were able to
choose fruitful methods and evade pitfalls \cite[p.~188]{Grab}
\end{quote}  
How can deep insight manage to ``evade pitfalls'' if the foundations
are regarded as inconsistent?  Grabiner further claims that,
\begin{quote}
[s]ince an adequate response to Berkeley's objections would have
involved recognizing that an equation involving limits is a shorthand
expression for a sequence of inequalities---a subtle and difficult
idea---no eighteenth century analyst gave a fully adequate answer to
Berkeley \cite[p.~189]{Grab}.
\end{quote}
This is an astonishing claim, which amounts to reading back into
history, developments that came much later.  Such a claim amounts to
postulating the inevitability of a triumphant march, from Berkeley
onward, toward the radiant future of Weierstrassian epsilontics.%
\footnote{Hourya Benis Sinaceur similarly mocks ``l'hypoth\`ese d'un
d\'eveloppement monolithique, non diff\'erenci\'e, continue, univoque
de l'analyse au XIX$^e$ siecle.  Le filon, perdu avec Bolzano, est
heureusement retrouv\'e par Cauchy, qui le red\'ecouvre, l'\'elargit,
l'exploite plus amplement, le transmet \`a Weierstrass, etc.''
\cite[p.~103]{Si}, and names Boyer \cite[p.~271]{Boy59} as one of the
culprits in perpetuating such a myth.}
The claim of such inevitability in our opinion is an assumption that
requires further argument.

Berkeley was, after all, attacking the coherence of {\em
infinitesimals\/}.%
\footnote{Berkeley's criticism is dissected into its metaphysical and
logical components by D.~Sherry \cite{She87}.}
He was {\em not\/} attacking the coherence of some kind of incipient
form of Weierstrassian epsilontics and its inequalities.  Isn't there
a simpler answer to Berkeley's query, in terms of a passage from a
point of B-continuum (see Appendix~\ref{rival}), to the infinitely
close point of the A-continuum, namely passing from a variable
quantity to its limiting constant quantity?

A related attitude on the part of Felscher is
discussed in Section~\ref{two}.  

Like Felscher, Grabiner \cite[p.~190]{Grab} suppresses Cauchy's
reference to an infinitesimal increment of the independent variable
when citing Cauchy's first definition (see Table~\ref{continuity} in
Section~\ref{two} above), thereby managing to avoid discussing
Cauchy's infinitesimals altogether.  We encounter the oft-repeated
claim about ``the same confusion between uniform and point-wise
convergence''%
\footnote{See Subsection~\ref{micro1} for a detailed discussion of the
controversy over the ``sum theorem''.}
\cite[p.~191]{Grab}. Her discussion of Cauchy's rigor does not mention
that what rigor meant to Cauchy was the replacement of the principle
of the generality of algebra, by geometry, including infinitesimals.
Grabiner correctly points out \cite[p.~193]{Grab} that
``Mathematicians are used to taking the rigorous foundations for
calculus for granted.''  She concludes: ``What I have tried to do as a
historian is to reveal what went into making up that great
achievement.''  What we have tried to do is to introduce a necessary
correction to a post-Weierstrassian reading of Cauchy, influenced by
an automated infinitesimal-to-limits translation originating no later
than Boyer \cite[p.~277]{Boy}, see \cite{KK11b} for additional
details.

\subsection{Devlin: 
Will the real Cauchy--Weierstrass please stand up?}

The following exchange between the second-named author and Keith
Devlin, co-founder and Executive Director of Stanford University's
H-STAR institute, took place in may 2011.  The exchange was occasioned
by Devlin's article ``Will the real continuous function please stand
up?''%
\footnote{The article is online at
http://www.maa.org/devlin/devlin\_11\_06.html}
The article refers to an alleged ``Cauchy--Weierstrass definition of
continuity''.

\bigskip\noindent MK: I read with interest your online article on
``continuous function please stand up''.  I am not sure which
Cauchy--Weierstrass definition you are referring to.  What exactly
{\em was\/} Cauchy's definition of continuity?  I think I know what
Weierstrass's was.

\medskip\noindent
KD: See%
\footnote{The link provided is a link to an electronic version of
Grabiner's article \cite{Grab}.}
http://www.maa.org/pubs/Calc\_articles/ma002.pdf

\medskip\noindent MK: But Grabiner does not say that Cauchy gave an
epsilontic definition of continuity. I wonder if there is a reason for
that.  Do you know what Cauchy's definition of continuity was?

\medskip\noindent KD: Actually, she does say that: ``Delta-epsilon
proofs are first found in the works of Augustin-Louis Cauchy
(1789-1867).  This is not always recognized, since Cauchy gave a
purely verbal definition of limit, which at first glance does not
resemble modern definitions."  In the UK when I was a student, the
standard term to refer to the usual~$\epsilon$-$\delta$ definition was
``the Cauchy--Weierstrass definition", and by my read of Grabiner that
description seems appropriate.

\medskip\noindent MK: But my question concerned the definition of {\em
continuity\/}.  In the UK, did they refer to a Cauchy--Weierstrass
definition of continuity?

\medskip\noindent
KD: It's the same.~$f(x)$ is continuous at~$a$ if and only if~$f(a)$
is defined and~$\lim_{x\to a}f(x) = f(a)$.

\medskip\noindent MK: Keith, this is {\em your\/} definition of
continuity.  It is {\em not\/} Cauchy's definition of continuity, for
several reasons: (1) Cauchy did not at any time work with ``continuity
at a point", it is always continuity ``between two limits", i.e. on an
interval. (2) Cauchy's definition is invariably an infinitesimal one:
``an infinitesimal~$x$-increment always produces an
infinitesimal~$y$-increment".  (3) Elsewhere, Cauchy is working for
the most part with a kinetic definition of limit, akin to that found
in Newton.%
\footnote{See footnote~\ref{pourciau}.}
Thus, the primary notion is that of a variable quantity.  Both
infinitesimals and limits are defined in terms of variable quantity.
Could the notion of a ``Cauchy--Weierstrass definition of continuity"
be an ahistorical blunder, regardless of whether it was adopted by the
mathematicians at a UK college?

\medskip\noindent
KD: That was actually the way continuity was taught in the UK.  To my
mind, once you have the basic epsilon-delta idea, which is a clever
static capture of the dynamic concept, everything just drops out. It
comes down to that one idea. And according to that Grabiner article,
that was due to Cauchy and then Weierstrass.  Of course, once we had
that definition, our entire conception changed, and it may be that it
is only from a modern perspective that we can clearly see that there
is really just one fundamental notion there. That happens all the
time, of course. Our present day conception of Newton's differential
calculus is not the one Newton had.%
\footnote{Pourciau refutes this view; see footnote~\ref{pourciau}.}
Sounds like you are trying to get at how the people at the time
conceived things.

\medskip\noindent MK: What you seem to be saying is that, even though
Cauchy himself said nothing of the sort regarding continuity, from a
certain ``modern perspective" is it proper to interpret Cauchy in a
certain way, now that in retrospect we know that the epsilontic idea
is the ``fundamental notion".  This is precisely the ideology that I
seek to refute.  Felix Klein pointed out in 1908 (in the original
edition of his ``Elementary mathematics from an advanced viewpoint")
that there are really {\em two\/} parallel strands to the development
of analysis, the epsilontic one and the infinitesimal one.  Emile
Borel in 1902 developed further du Bois-Reymond's theory of
infinitesimal-enriched continua in terms of growth rates of functions,
and specifically refers to Cauchy's work in this direction from 1829
as inspiration.  Note that Klein and Borel said this well before
Robinson.  Viewed as the ``father" of the infinitesimal approach
ultimately vindicated by Robinson, Cauchy needn't suffer a forced
explanation as a proto-Weierstrassian.  Rather, his approach to
infinitesimals was further developed by Stolz, du Bois-Reymond, Borel,
Levi-Civita, and others around the turn of the century.  When Skolem
developed first non-standard models of arithmetic in 1934, he was
inspired partly by du Bois-Reymond's work, according to Robinson.  In
1948, Hewitt first developed hyperreal fields using ultraproducts.  In
1955, \Los{} proved his theorem for ultraproducts, which implies the
transfer principle, which is the mathematical implementation of
Leibniz's heuristic ``law of continuity": whatever succeeds for the
finite, should also succeed for the infinite.  Isn't this lineage of
infinitesimal calculus more plausible than an alleged
``Cauchy--Weierstrass definition of continuity"?

\medskip\noindent
KD: Ah, I assume you are a historian of mathematics. For sure my
perspective is very much that of mathematicians, who if they have any
interest in history (and I do) it is precisely from the perspective
you describe, tracing back the origins of the now accepted concepts in
terms that we now use.  For instance, we interpret Newton's work on
calculus as applied to continuous functions of a real variable, when
that is not at all what Newton was doing. (Actually, I don't think we
can really understand what he was doing, since we are rooted in a
modern perspective.) Good luck with your investigations.%
\footnote{Devlin's article acknowledges that it is based on N\'u\~nez
{\em et al\/} \cite{NEM}, who similarly employ the term
``Cauchy--Weierstrass definition of continuity''.}

\bigskip

To put the above exchange in perspective, it may be useful to recall
Grattan-Guinness's articulation of a historical reconstruction project
in the name of H.~Freudenthal~\cite{Fr71b}, in the following terms:
\begin{quote}
it is mere feedback-style ahistory to read Cauchy (and contemporaries
such as Bernard Bolzano) as if they had read Weierstrass already.  On
the contrary, their own pre-Weierstrassian muddles%
\footnote{\label{muddle}Grattan-Guinness's term ``muddle'' refers to
an irreducible ambiguity of historical mathematics such as Cauchy's
sum theorem of 1821.}
need historical reconstruction \cite[p.~176]{Grat04}.
\end{quote}

\subsection{Grattan-Guinness on beautiful mathematics}
\label{gg2}

On the subject of Robinson's theory, Grattan-Guinness comments as
follows:
\begin{quote}
I made no mention of non-standard analysis in my book, for it was
obvious to me that this very beautiful piece of mathematics had
nothing to tell us historically \cite[p.~247]{Grat78}.
\end{quote}
When Grattan-Guinness announced that Robinson's construction of a
non-Archimedean extension of the reals ``bears no resemblance to past
arguments in favor of infinitesimals'' \cite[p.~247]{Grat78}, he was
only telling part of the story.  True, the construction favored by
Robinson exploited powerful compactness theorems%
\footnote{Robinson uses instead the term ``finiteness principle of
lower predicate calculus'' \cite[p.~48]{Ro66} for what is known today
as the compactness theorem.}
originating with Malcev~\cite{Ma} and eschewed the sequential
approach.  On the other hand, the ultrapower construction of the
hyperreals, pioneered by E.~Hewitt \cite{Hew48} in 1948 and
popularized by Luxemburg \cite{Lu62} in 1962, is firmly rooted in the
sequential approach, and hence connects well with the kinetic vision
of Cauchy, shared by L.~Carnot.%
\footnote{Some details on the ultrapower construction appear in
Appendix~\ref{rival}.}

\section{Timeline of modern infinitesimals from Cauchy onward}
\label{time}

The historical sequence of events, as far as continuity is concerned,
was as follows:%
\footnote{We leave out Bolzano's contribution which, while prior to
Cauchy's, did not exert any influence until the 1860s.}
\begin{enumerate}
\item
first came Cauchy's infinitesimal definition, namely ``infinitesimal
$x$-increment always produces an infinitesimal~$y$-increment";
\item
then came the Dirichlet/Weierstrass-style nominalistic reconstruction
of the original definition in terms of real inequalities, dispensing
with infinitesimals;
\item
Cauchy's work influenced investigations in infinitesimal-enriched
continua by du Bois-Reymond, E.~Borel, and others at the turn of the
century, see Table~\ref{heuristic} (a more detailed account appears
below).
\end{enumerate}

The historical {\em priority\/} of the infinitesimal definition is
clear; what is open to debate is the role of the modern definition in
{\em interpreting\/} the historical definition.  The common element
here is the {\em null sequence\/}, a basis both for Cauchy
infinitesimals, and, via intermediate developments in growth rates of
functions as explained below, for ultrapower-based infinitesimals.%
\footnote{\label{martin}Modulo suitable foundational material, one can
ensure that every hyperreal infinitesimal is represented by a null
sequence; an appropriate ultrafilter (called a P-point) will exist if
one assumes the continuum hypothesis, or even the weaker Martin's
axiom (see Cutland {\em et al\/} \cite{CKKR} for details).}

\renewcommand{\arraystretch}{1.3}
\begin{table}
\[
\begin{tabular}[t]
{ | p{.4in} || p{1.4in} | p{2.7in} | } \hline years & author &
contribution \\ \hline\hline 1821 & Cauchy & Infinitesimal definition
of continuity \\ \hline 1827 & Cauchy & Infinitesimal delta function
\\ \hline 1829 & Cauchy & Defined ``order of infinitesimal'' in terms
of rates of growth of functions \\ \hline 1852 & Bj\"orling & Dealt
with convergence at points ``indefinitely close'' to the limit \\
\hline 1853 & Cauchy & Clarified hypothesis of ``sum theorem'' by
requiring convergence at infinitesimal points \\ \hline 1870-1900 &
Stolz, du~Bois-Reymond, and others & Infinitesimal-enriched number
systems defined in terms of rates of growth of functions\\ \hline 1902
& Emile Borel & Elaboration of du Bois-Reymond's system \\ \hline 1910
& G.~H.~Hardy & Provided a firm foundation for du Bois-Reymond's
orders of infinity \\ \hline 1926 & Artin--Schreier & Theory of real
closed fields \\ \hline 1930 & Tarski & Existence of ultrafilters \\
\hline 1934 & Skolem & Nonstandard model of arithmetic \\ \hline 1948
& Edwin Hewitt & Ultrapower construction of hyperreals \\ \hline 1955
& \Los{} & Proved \Los's theorem forshadowing the transfer principle
\\ \hline 1961, 1966 & Abraham Robinson & Non-Standard Analysis \\
\hline 1977 & Edward Nelson & Internal Set Theory \\ \hline
\end{tabular}
\]
\caption{\textsf{Timeline of modern infinitesimals from Cauchy to
Nelson}.}
\label{heuristic}
\end{table}
\renewcommand{\arraystretch}{1}

\subsection{From Cauchy to du Bois-Reymond}

Cauchy's theory of arbitrary orders of magnitude for his
infinitesimals was a harbinger, not of Weierstrassian epsilontics, but
of later theories of infinitesimal-enriched continua as developed by
Stolz \cite{Stolz}, du Bois-Reymond \cite{Boi}, and others.  Emile
Borel's appreciation of the essential continuity between Cauchy's
theory and that of his late 19th century heirs was already discussed
in Subsection~\ref{two}.%
\footnote{See the main text following footnote~\ref{ehr}.}
Robinson \cite[p.~277-278]{Ro66} traces the evolution of the
infinitesimal ideal from Cauchy to du Bois-Reymond and Stolz, to
Skolem.%
\footnote{\label{1966} See the main text around footnote~\ref{1966a}.}

\subsection{From du Bois-Reymond to Robinson}

Du Bois-Reymond's investigations were in turn pursued further by such
mathematicians as Emile Borel.  In 1902, Borel \cite[p.~35-36]{Bo02}
cites Cauchy's definition of such ``order of infinitesimal"%
\footnote{This definition appears in Oeuvres de Cauchy, s\'erie 2,
tome 4, p.~181, corresponding to Cauchy's {\em Le\c cons sur le calcul
diff\'erentiel\/} from 1829, appearing in a section entitled
``Pr\'eliminaires"; see G.~Fisher \cite[p.~144]{Fi81}.}
as inspiration for du Bois-Reymond's theory.  Hardy \cite{Har}
provided a firm foundation for du Bois-Reymond's orders of infinity in
1910.  Artin and Schreier \cite{AS} developed the theory of real
closed fields in 1926.  Skolem \cite{Sk} constructed the first
non-standard models of arithmetic in 1934.  Hewitt's hyperreals of
1948 and \Los's theorem of 1955 have already been discussed.  In the
1960s, Robinson proposed an infinitesimal theory \cite{Ro66} as an
alternative to Weierstrassian epsilontic analysis.  The evolution of
infinitesimal from Cauchy to Robinson is summarized in
Table~\ref{heuristic}.

Edward Nelson \cite{Ne} in 1977 proposed an axiomatic theory parallel
to Robinson's theory, see Subsection~\ref{nelson}.  P.~Ehrlich
recently constructed an isomorphism of maximal hyperreal and surreal
fields, resulting in a ``unification of all numbers great and small''
\cite{Eh12}.

\subsection{Nelson}
\label{nelson}

Edward Nelson \cite{Ne} in 1977 proposed an axiomatic theory parallel
to Robinson's theory.  Nelson's axiomatisation, called Internal Set
Theory (IST), takes the form of an enrichment of the Zermelo-Fraenkel
set theory (ZFC) \cite{Fran}.  IST amounts to a more stratified
axiomatisation for set theory, more congenial to infinitesimals.

In Nelson's system, the usual construction of~$\R$, when interpreted
with respect to the foundational background provided by IST, produces
a number system already possessing entities behaving as
infinitesimals.

In more detail, Nelson's approach is a re-thinking of the foundational
material with a view to allowing a more stratified (hierarchical)
number line.  Thus, the canonical set theory, namely ZFC, is modified
by the introduction of a unary predicate ``standard''.  Then what is
known as the usual construction of the ``real'' line produces a line
that bears a striking resemblance to the Hewitt-\L o\'s-Robinson
hyperreals.

To illustrate the power of Nelson's approach, we quote Alain Robert
who starts his book on IST \cite{Robert03} with a homage to Leonard
Euler:

\begin{quotation}
\small\noindent Here is how [Euler] deduces the expansion of the
cosine function [\dots ] He starts from the de Moivre
formula\begin{eqnarray*} \cos nz &=& \frac{1}{2}\left[(\cos z + i\sin
z)^n + (\cos z - i\sin z)^n\right]\\ &=& \cos^n z
-\frac{n(n-1)}{1\cdot 2} \cos^{n-2}z \sin^2 z\\ &&\quad +
\frac{n(n-1)(n-2)(n-3)}{1\cdot 2 \cdot 3\cdot 4} \cos^{n-4}z \sin^4 z+
\cdots \end{eqnarray*} and writes [\dots ]
\begin{quotation}\noindent \emph{sit arcus~$z$ infinite parvus; erit
$\cos\cdot z =1$,~$\sin\cdot z = z$; sit autem~$n$ numerus infinite
magnus, ut sit arcus~$nz$ finitae magnitudinis, puta}~$nz=v$
\[
\cos \cdot v = 1 -\frac{v^2}{2!} + \frac{v^4}{4!}\mbox{---etc.}
\]
\end{quotation} 
\normalsize
\end{quotation} 
In the context of internal set theory, this argument becomes fully
rigorous.

\section{Conclusion}
\label{five}

As we have seen, Cauchy experimented with a range of foundational
approaches, including infinitesimal methodologies.  He anticipated a
number of mathematical developments that occurred decades later.  One
such development was the study of the asymptotic behavior and rates of
growth of functions, culminating in the construction of
infinitesimal-enriched continua by du Bois-Reymond, Emile Borel,
G.~H.~Hardy, and others.  Another development anticipated by Cauchy
was the Dirac delta function, of which Cauchy gave an infinitesimal
definition.  Similarly, Bj\"orling and Cauchy both articulated
single-variable definitions of stronger forms of continuity and
convergence (today called {\em uniform\/}), in the context of
infinitesimal-enriched continua.  One avenue for further exploration
is Cauchy's infinitesimal approach to the degrees of contact among
curves.%
\footnote{Cauchy does not hesitate to characterize the center of
curvature as the meeting point of a pair of infinitely close normals
\cite[p. 91]{Ca26}.}
Such visionary anticipations in Cauchy are studiously obfuscated by
scholars sporting protective blinders, conditioned by a conceptual
framework stemming from a nominalistic reconstruction of analysis set
in motion by the ``great triumvirate'' of Cantor, Dedekind, and
Weierstrass.  According to the rules of such feedback-style ahistory,
to borrow Grattan-Guinness's term,%
\footnote{See main text at footnote~\ref{muddle}.}
Cauchy's infinitesimals are subjected to an automated translation to
limits, the latter being promoted to first fiddle, Cauchy's explicit
statements to the contrary notwithstanding (see
Subsection~\ref{gilain}).  Yet, many historians have refused to toe
the triumvirate line.  Cauchy's anticipations are highlighted in
insightful studies by Freudenthal \cite{Fr71a}, Robinson~\cite{Ro66},
Hourya Benis Sinaceur \cite{Si}, Lakatos \cite{La78}, Cleave
\cite{Cl79}, Cutland {\em et al\/} \cite{CKKR}, Medvedev \cite{Me87},
Laugwitz \cite{Lau89}, Sad {\em et al\/} \cite{STB}, Br\aa ting
\cite{Br07}, and others.

It bears pointing out that the the calculus of Newton and Leibniz was
not the same as ours, since they did not have a continuum that lives
up to modern standards.  Yet historians routinely attribute the
invention of the calculus to them, implicitly taking it for granted
that the calculus of Newton and Leibniz can be interpreted in terms of
modern calculus, and vice versa.  Isn't it time we applied such
bi-interpretability to Cauchy's infinitesimals, as well?

\section*{Acknowledgments}

We are grateful to Hourya Benis Sinaceur, David Sherry, and Detlef
Spalt for a careful reading of an earlier version of the manuscript
and for their helpful comments.

\appendix

\section{Rival continua}
\label{rival}

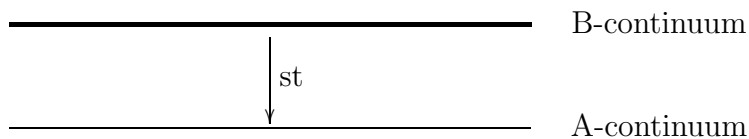
\begin{figure}
\[
\xymatrix@C=95pt{{} \ar@{-}[rr] \ar@{-}@<-0.5pt>[rr]
\ar@{-}@<0.5pt>[rr] & {} \ar@{->}[d]^{\hbox{st}} & \hbox{\quad
B-continuum} \\ {} \ar@{-}[rr] & {} & \hbox{\quad A-continuum} }
\]
\caption{Thick-to-thin: taking standard part (the thickness of the top
line is merely conventional)}
\label{31}
\end{figure}

A Leibnizian definition of the derivative as the infinitesimal ratio
\[
\frac{\Delta y}{\Delta x},
\] 
whose logical weakness was criticized by Berkeley, was modified by
A.~Robinson by exploiting a map called {\em the standard part\/},
denoted~``st'', from the finite part of a B-continuum (for
``Bernoullian''), to the A-continuum (for ``Archimedean''), as
illustrated in Figure~\ref{31}.%
\footnote{In the context of the hyperreal extension of the real
numbers, the map ``st'' sends each finite point~$x$ to the real point
st$(x)\in \R$ infinitely close to~$x$.  In other words, the map ``st''
collapses the cluster of points infinitely close to a real number~$x$,
back to~$x$.  A comparative study of continua from a predicative angle
may be found in Feferman~\cite{Fef}.}

We will denote such a B-continuum by a new symbol \RRRhuge.  We will
also denote its finite part, by
\[
\RRR_{<\infty} = \left\{ x\in \RRR : \; |x|<\infty \right\};
\]
namely, the difference $\RRR \setminus \RRR_{<\infty}$ consists of the
inverses of nonzero infinitesimals.  The map ``st'' sends each finite
point~$x\in \RRR$, to the real point st$(x)\in \R$ infinitely close
to~$x$:
\[
\xymatrix{\quad \RRR_{{<\infty}}^{~} \ar[d]^{{\rm st}} \\ \R}
\]

Robinson's answer to Berkeley's {\em logical criticism\/} (see
D.~Sherry \cite{She87}) is to define the derivative of~$y=f(x)$ as
\[
f'(x)= \text{st} \left( \frac{\Delta y}{\Delta x} \right),
\]
rather than the infinitesimal ratio~$\Delta y/\Delta x$ itself, as in
Leibniz. ``However, this is a small price to pay for the removal of an
inconsistency'' \cite[p~266]{Ro66}.

\begin{figure}
\includegraphics[height=2in]{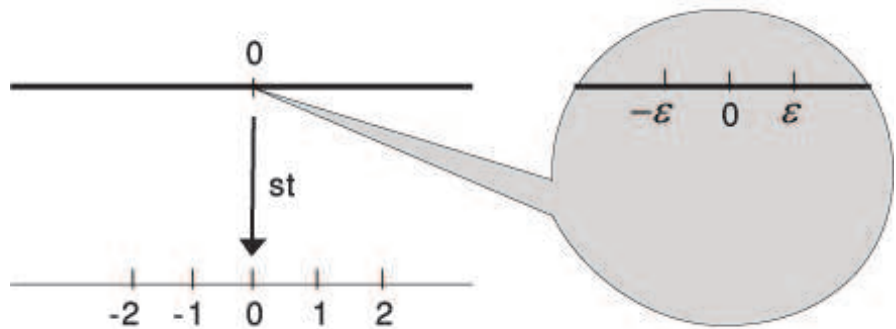}
\caption{Zooming in on infinitesimal~$\epsilon$}
\label{tamar}
\end{figure}

We illustrate the construction by means of an infinite-resolution
microscope in Figure~\ref{tamar}.

Note that both the term ``hyper-real field'', and an ultrapower
construction thereof, are due to E.~Hewitt in 1948, see
\cite[p.~74]{Hew48}.  The transfer principle allowing one to extend
every first-order real statement to the hyperreals, is due to
J.~{\L}o{\'s} in 1955, see \cite{Lo}.  Thus, the Hewitt-{\L}o{\'s}
framework allows one to work in a B-continuum satisfying the transfer
principle.  

Hewitt's construction of hyper-real fields has roots in functional
analysis, including works by Gelfand and Kolmogorov \cite{GK}.  In
1990, Hewitt reminisced about his
\begin{quote}
``efforts to understand the ring of all real-valued continuous [not
necessarily bounded] functions on a completely regular~$T_0$-space. I
was guided in part by a casual remark made by Gel'fand and Kolmogorov%
\footnote{Here Hewitt cites Gelfand and Kolmogoroff \cite{GK}.}
[\dots]  Along the way I found a novel class of real-closed fields that
superficially resemble the real number field and have since become the
building blocks for nonstandard analysis. I had no luck in talking to
Artin about these hyperreal fields, though he had done interesting
work on real-closed fields in the 1920s. (My published `proof' that
hyperreal fields are real-closed is false: John Isbell earned my
gratitude by giving a correct proof some years later.)  [\dots]  My
ultra-filters also struck no responsive chords. Only Irving Kaplansky
seemed to think my ideas had merit. My first paper on the subject was
published only in 1948'' \cite{He90}
\end{quote}
(Hewitt goes on to detail the influence of his 1948 text).  Here
Hewitt is referring to Isbell's 1954 paper \cite{Isb}, proving that
Hewitt's hyper-real fields are real closed.  Note that a year later,
\Los~\cite{Lo} proved the general transfer principle for such fields,
implying in particular the property of being real closed.

To elaborate on the ultrapower construction of the hyperreals,
let~$\Q^\N$ denote the space of sequences of rational numbers.
Let~$\left( \Q^\N \right)_C$ denote the subspace consisting of Cauchy
sequences.  The reals are by definition the quotient field
\[
\R:= \left. \left( \Q^\N \right)_C \right/ \mathcal{F}_{\!n\!u\!l\!l},
\]
where the ideal~$\mathcal{F}_{\!n\!u\!l\!l}$ contains all the null
sequences.  Meanwhile, an infinitesimal-enriched field extension
of~$\Q$ may be obtained by forming the quotient
\[
\left.  \Q^\N \right/ \mathcal{F}_{u},
\]
see Figure~\ref{helpful}.  Here a sequence~$\langle u_n : n\in \N
\rangle$ is in~$\mathcal{F}_{u}$ if and only if the set
\[
\{ n \in \N : u_n = 0 \}
\]
is a member of a fixed ultrafilter.%
\footnote{An ultrafilter on~$\N$ can be thought of as a way of making
a systematic choice, between each pair of complementary infinite
subsets of~$\N$, so as to prescribe which one is ``dominant'' and
which one is ``negligible''.  Such choices have to be made in a
coherent manner, e.g., if a subset~$A\subset \N$ is negligible then
any subset of~$A$ is negligible, as well.  The existence of
ultrafilters was proved by Tarski \cite{Tar}, see Keisler
\cite[Theorem~2.2]{Ke08}.  See a related remark about P-points in
footnote~\ref{martin}.}
To give an example, the sequence~$\left\langle \tfrac{(-1)^n}{n}
\right\rangle$ represents a nonzero infinitesimal, whose sign depends
on whether or not the set~$2\N$ is a member of the ultrafilter.  To
obtain a full hyperreal field, we replace~$\Q$ by~$\R$ in the
construction, and form a similar quotient
\[
\RRR:= \left.  \R^\N \right/ \mathcal{F}_{u}.
\]

\begin{figure}
\[
\xymatrix{ && \left( \left. \Q^\N \right/ \mathcal{F}_{\!u}
\right)_{<\infty} \ar@{^{(}->} [rr]^{} \ar@{->>}[d]^{\rm st} &&
\RRR_{<\infty} \ar@{->>}[d]^{\rm st} \\ \Q \ar[rr] \ar@{^{(}->} [urr]
&& \R \ar[rr]^{\simeq} && \R }
\]
\caption{\textsf{An intermediate field~$\left. \Q^\N \right/
\mathcal{F}_{\!u}$ is built directly out of~$\Q$}}
\label{helpful}
\end{figure}

A more detailed discussion of the ultrapower construction can be found
in M.~Davis~\cite{Da77}.  See also B\l aszczyk \cite{Bl} for some
philosophical implications.  More advanced properties of the
hyperreals such as saturation were proved later, see Keisler
\cite{Kei} for a historical outline.  A helpful ``semicolon'' notation
for presenting an extended decimal expansion of a hyperreal was
described by A.~H.~Lightstone~\cite{Li}.  See also P.~Roquette
\cite{Roq} for infinitesimal reminiscences.  A discussion of
infinitesimal optics is in K.~Stroyan \cite{Str},
H.~J.~Keisler~\cite{Ke}, D.~Tall~\cite{Ta80}, L.~Magnani \&
R.~Dossena~\cite{MD, DM}, and Bair \& Henry \cite{BH}.

Applications of the B-continuum range from aid in teaching calculus
\cite{El, KK1, KK2, Ta91, Ta09a} to the Bolzmann equation (see
L.~Arkeryd~\cite{Ar81, Ar05}); modeling of timed systems in computer
science (see H.~Rust \cite{Rust}); mathematical economics (see
R.~Anderson \cite{An00}); mathematical physics (see Albeverio {\em et
al.\/} \cite{Alb}); etc.

An interesting elementary example of a B-continuum is the ring of dual
numbers, that is, numbers of the form~$a + b\delta$,
where~$\delta^2=0$.  It is useful, for example, in the theory of
linear algebraic groups and linear Lie groups, and enable a quick and
transparent computation of Lie algebras of groups
like~$\text{SO}_3(\mathbb{R})$.

Recently, Giordano \cite{Gio10a, Gio10b, Gio11a} introduced a
systematic way of enriching the reals by nilpotent infinitesimals,
referring to the resulting structure as \emph{the ring of Fermat
reals\/}.

The transfer principle of the modern theory of infinitesimals is a
mathematical implementation of Leibniz's heuristic law of continuity,
see Robinson \cite[p.~266]{Ro66}, and Laugwitz \cite{Lau92}.  The
standard part function is a mathematical implementation of Fermat's
concept of adequality, see Table~\ref{heuristic2}.

\renewcommand{\arraystretch}{1.3}
\begin{table}
\[
\begin{tabular}[t]
{ | p{2.3in} || p{2.3in} | } \hline heuristic concept & mathematical
implementation \\ \hline\hline adequality & standard part function \\
\hline law of continuity & transfer principle \\ \hline
infinitesimal-enriched continuum & hyperreal number line \\ \hline
\end{tabular}
\]
\caption{\textsf{Modern mathematical implementation of 17th century
heuristic concepts}.}
\label{heuristic2}
\end{table}
\renewcommand{\arraystretch}{1}


\bigskip\noindent \textbf{Alexandre Borovik} is Professor of Pure
Mathematics at the University of Manchester, United Kingdom, where he
has been working for the past 20 years.  His principal research lies
in algebra, model theory, and combinatorics---topics on which he
published several monographs and a number of papers.  Recently he has
become interested in studying algebraic phenomena inhabiting the murky
boundary between finite and infinite.  He also has an interest in
cognitive aspects of mathematical practice and recently published a
book \emph{Mathematics under the Microscope: Notes on Cognitive
Aspects of Mathematical Practice} which explains a mathematician's
outlook at psychophysiological and cognitive issues in mathematics. An
(almost final) draft of his book \emph{Shadows of the Truth:
Metamathematics of Elementary Mathematics} can be found at
http://www.maths.manchester.ac.uk/$\sim$avb/ST.pdf

\medskip\noindent\textbf{Mikhail G. Katz} is Professor of Mathematics
at Bar Ilan University, Ramat Gan, Israel.  Two of his joint studies
with Karin Katz were published in {\em Foundations of Science\/}: ``A
Burgessian critique of nominalistic tendencies in contemporary
mathematics and its historiography" and ``Stevin numbers and reality",
online respectively at

http://dx.doi.org/10.1007/s10699-011-9223-1 and at

http://dx.doi.org/10.1007/s10699-011-9228-9

A joint study with Karin Katz entitled ``Meaning in classical
mathematics: is it at odds with Intuitionism?" is due to appear in
{\em Intellectica\/}.

His joint study with David Tall, entitled ``The tension between
intuitive infinitesimals and formal mathematical analysis", is due to
appear as a chapter in a book edited by Bharath Sriraman, see

\noindent
http://www.infoagepub.com/products/Crossroads-in-the-History-of-Mathematics

\end{document}